\documentstyle[12pt]{article}
\begin{document}
%\large
\title{Unbounded operators on Banach spaces over the quaternion field}
\author{
%\par \vspace{0.5cm}
%\par {{\large
S.V. Ludkovsky}
\date{27.08.2003}
\maketitle
\par \vspace{0.5cm}
\footnote{This work was partially supported by DAAD, Germany.}
\par \vspace{0.5cm}
\section{Introduction}
The quaternion field is the algebra over $\bf R$, but it is not the algebra
over $\bf C$, since each embedding of $\bf C$ into $\bf H$
is not central. Therefore, the investigation of operator algebras over
$\bf H$ can not be reduced to algebras of operators over $\bf C$.
On the other hand, the developed below theory of operator algebras over
$\bf H$ has many specific features in comparison with the general
theory of operator algebras over $\bf R$ due to the graded structure
of $\bf H$. Results of this work may be also used for the development
of non-commutative geometry, superanalysis, quantum mechanics
over $\bf H$ and the representation theory of non locally compact
groups such as groups of diffeomorphisms and loops of quaternion manifolds
(see \cite{connes,oystaey,emch,lulgcm,lupm}).
The vast majority of previous works on superanalysis was devoted
to supercommutative superalgebras of the type of Grassman algebra,
but for the non-commutative superalgebras it was almost undeveloped.
The quaternion field serves as very important example of the superalgebra,
which is not supercommutative.
In this work the results of previous works on this theme of the author
are used, in particular, non-commutative integral over $\bf H$
\cite{luoyst,lufsqv} that serves the analog of the Cauchy-type integral
known for $\bf C$. Examples of quaternion unbounded operators are
differential operators and among them in partial derivatives.
They arise in the natural way, for example, the Klein-Gordon-Fock
equation can be written in the form $(\partial ^2/
\partial z^2+\partial ^2/ \partial {\tilde z}^2)f=0$
on the space of quaternion locally $(z,{\tilde z})$-analytical functions
$f$, where $z$ - is the quaternion variable, $\tilde z$ - is the adjoint
variable, $z{\tilde z}=|z|^2$. The Dirack operator for spin systems
over $\bf H^2$ can be written in the form
${ {0\quad \partial / \partial z}
\choose {{- \partial / \partial {\tilde z}} \quad 0}} $, that is used
in the theory of spin manifolds \cite{lawmich}, but each spin
manifold can be embedded into the quaternion manifold \cite{lufsqv}.
In this article main features of the quaternion case are given, since
in this article it is impossible to give the same broad theory
over $\bf H$, as well-developed theory of operators over $\bf C$
\cite{danschw,kadring}.
\section{Theory of unbounded operators}
{\bf 2.1. Definitions and Notes.}
Let $X$ be a Banach space (BS) over the quaternion field $\bf H$,
that is, $X$ is the additive group, multiplications of vectors $v\in X$
on scalars $a, b\in \bf H$ on the left and right satisfy axioms
of associativity and distributivity, there exists the norm
$\| v \|_X=:\| v \| $ on $X$ relative to which, $X$
is complete, where $\| av \| =|a|_{\bf H}\| v \| $,
$\| vb \| =|b|_{\bf H}\| v \| $, $\| v+w\| \le \| v \| +\| w \| $
for each $v, w\in X$, $a, b\in \bf H$.
Then $X$ has also the structure $X_{\bf R}$ of BS over $\bf R$, since
$\bf H$ is the algebra over $\bf R$ of dimension $4$.
\par An operator $T$ on a dense vector subspace ${\cal D}(T)$
in $X$ with values in BS $Y$ over $\bf H$ is called ($\bf H$-)right-linear
(RLO), if $(Tva)b=T(vab)$, $T(va+wb)=(Tv)a+(Tw)b$ for each $a, b\in \bf H$
and each $v, w\in {\cal D}(T)$ and in addition $T$ is $\bf R$-linear
on ${\cal D}(T)_{\bf R}$.
FOR RLO we also write $Tv$ instead of $T(v)$.
If $T$ is $\bf R$-linear and $b(Tav)=T(bav)$, $T(av+bw)=
a(Tv)+b(Tw)$ for each $a, b\in \bf H$, then $T$ is called
($\bf H$-)left-linear (LLO).
For LLO we also write $vT$ instead of $T(v)$.
An operator $T$ is called ($\bf H$-)linear, if it is LLO and RLO
simultaneously. An operator $T: {\cal D}(T)\to {\cal R}(T)$ we call
($\bf H$-)quasi-linear (QLO), if it is additive $T(v+w)=T(v)+T(w)$
and $\bf R$-homogeneous $T(av)=aT(v)$ for each $a\in \bf R$, $v$ and
$w\in X$, where ${\cal R}(T)\subset Y$ denotes the range of values
of this operator.
For example, products of quaternion holomorphic functions
are $\bf H$-quasi-linear (see \cite{luoyst}).
\par Let $L_q(X,Y)$ ($L_r(X,Y)$; $L_l(X,Y)$) denote BS of all
bounded QLO $T$ from $X$ into $Y$ (RLO aand LLO respectively),
$\| T\| :=\sup_{v\ne 0} \| Tv \| /\| v \| $; $L_q(X)
:=L_q(X,X)$, $L_r(X):=L_r(X,X)$ and $L_l(X):=L_l(X,X)$.
The resolvent set $\rho (T)$ of QLO $T$ is defined as
$\rho (T):=\{ z\in {\bf H}: \mbox{ there exists }
R(z;T) \in L_q(X) \} $, where $R(z;T):=R_z(T):=(zI-T)^{-1}$,
$I$ is the identity operator on $X$, $Iv=v$ for each $v\in X$,
analogously for RLO and LLO.
A spectrum is defined as $\sigma (T):={\bf H}\setminus \rho (T) $.
\par {\bf 2.2. Lemma.} {\it For each $z_1$ and $z_2\in \rho (T)$: \\
$(i)\quad R(z_2;T)-R(z_1;T)=R(z_1;T)(z_1-z_2)R(z_2;T)$.}
\par {\bf 2.3. Lemma.} {\it If $T$ is a closed QLO, then $\rho (T)$
is open in $\bf H$ and $R(z;T)$ is quaternion holomorphic on $\rho (T)$.}
\par {\bf Proof.} Let $z_0\in \rho (T)$, then
the operator $R(z_0;T)$ is closed and by the closed graph theorem
it is bounded, since it is defined everywhere. If
$z\in \bf H$ and $|z_0-z|< \| (z_0I-T)^{-1} \|^{-1}$, then
$(zI-T)=(z_0I-T)(I-R(z_0;T)(z_0-z))$, such that \\
$(i)\quad R(z;T)=
\{ \sum_{n=0}^{\infty }[R(z_0;T)(z_0-z)]^n \} R(z_0;T)\in L_s(X),$
since this series $[R(z_0;T)(z_0-z)]^n R(z_0;T)$ converges relative
to the norm topology in $L_s(X)$. In view of $(i)$ $R(z;T)$
is quaternion locally $z$-analytic, hence it is holomorphic on
$\rho (T)$ due to Theorem 2.16 \cite{luoyst}.
\par {\bf 2.4. Notes and Definitions.} For BS $X$ over
$\bf H$ its right-adjoint space $X^*_r$
is defined as consisting of all functionals $f: X\to \bf H$
such that $f$ is $\bf R$-linear and $\bf H$-right-linear.
Analogously we put $X^*_q:=L_q(X,{\bf H})$ and $X^*_l:=L_l(X,{\bf H})$,
where $X^*_q$ is the topologically quasi-adjoint space,
$X^*_l$ is the topologically left-adjoint space.
Then $X^*_s$ is BS over $\bf H$ with the norm
$\| f \| :=\sup_{x\ne 0}|fx|/ \| x \|$.
If $X$ and $Y$ are BS over $\bf H$ and $T: X\to Y$ belongs to $L_s(X,Y)$,
then $T^*$ is defined by the equation: $(T^*y^*)(x):=y^*\circ T(x)$
for each $y^*\in Y^*_s$, $x\in X$.
Then $T^*\in L_l(Y^*_r,X^*_r)$ for each $T\in L_r(X,Y)$.
If $T\in L_l(X,Y)$, then $T^*\in L_r(Y^*_l,X^*_l)$, since $y^*\circ T(x)=
xT^*y^*$ in the symmetrical notation, where $x\in X$.
If $T\in L_q(X,Y)$, then $T^*\in L_q(Y^*_q,X^*_q)$.
\par Let $\hat X$ and $\hat Y$ be images relative of the natural embedding
$X$ and $Y$ into $X^{**}$ and $Y^{**}$ respectively.
For each $T\in L_s(X,Y)$ we define ${\hat T}\in L_s({\hat X},{\hat Y})$
by the equation ${\hat T}({\hat x})=\hat y$, where $y=T(x)$.
For each function $S$ defined on the region $X^{**}\supset dom (S)\supset
\hat X$ and such that $S({\hat x})=T({\hat x})$ for each
${\hat x}\in \hat X$ is called the extension of $T$.
\par Lemmas $II.3.12; VI.2.2-4,6,7$ and corollary $II.3.13$
from \cite{danschw} are analogous in the case $\bf H$
instead of $\bf C$, taking in the proof of Lemma $II.3.12$
$\| z \| =\| y+\alpha x \| =|\alpha | \| \alpha ^{-1}y+x \|
\ge |\alpha | d$.
Take by induction ${\cal D}(T^n):= \{ x: x\in {\cal D}(T^{n-1}),
T^{n-1}(x)\in {\cal D}(T) \} $,
${\cal D}(T^{\infty }):=\bigcap_{n=1}^{\infty }{\cal D}(T^n)$,
where $T^0:=I$, $T^n(x):=T(T^{n-1}(x))$.
\par {\bf 2.5. Lemma.} {\it Let $T\in L_s(X)$, then
$\sigma (T^*)=(\sigma (T))^{\tilde .}$ and $(R(\lambda ,T))^*=
R(\lambda ^*,T^*)$, where $\lambda ^*I:=(\lambda I)^*$,
$s\in \{ q, r, l \} $.}
\par {\bf Proof.} If $S\in L_s(X,Y)$ and there exists
$S^{-1}\in L_s(Y,X)$, then $S^* \in L_u(Y^*_s,X^*_s)$ has the inverse
$(S^*)^{-1} \in L_u(X^*_s,Y^*_s)$ and $(S^{-1})^*=
(S^*)^{-1}$, where $(s,u)\in \{ (q,q); (r,l); (l,r) \} $.
Then $(\lambda I-T)^*y^*=y^*\circ (\lambda I-T)$
$=y^*\circ (\lambda I)-y^*\circ T$, consequently,
$(\lambda ^*I-T^*)[(\lambda I-T)^*]^{-1}=I$
and $(R(\lambda ,T))^* =R(\lambda ^*,T^*)$. 
\par {\bf 2.6. Definition and Note.} Denote by ${\cal H}(T)$
the family of all quaternion holomorphic functions (QHF) $f$ on
neighbourhoods $V_f$ for $\sigma (T)$, where $T\in L_s(X)$,
$s\in \{ q, r, l \} $, and for a QLO $T$ ${\cal H}_{\infty }(T)$
is a set of all QHF on neighbourhoods $U_f$ of $\sigma (T)$ and
$\infty $ in the one-point compactification $\bf {\hat H}$
of the quaternion field.
We choose a marked point $z_0\in \sigma (T)$.
For each $M=wJ+xK+yL\in \bf H_i$ with $|M|=1$, where $w, x, y\in \bf R$,
there exists a closed rectifiable path $\eta $ consisting of a finite
union of arches $\eta (s)=z_0+r_p\exp (2\pi sM)$ with
$s\in [a_p,b_p]\subset [0,1]\subset \bf R$ and segments of straight lines
$\{ z\in {\bf H}: z=z_0+(r_pt+r_{p+1}(1-t)) \exp (2\pi b_pM),
t\in [0,1] \} $ joining them, moreover $\eta \subset U\setminus \sigma (T)$,
where $a_p<b_p$ and $0<r_p<\infty $ for each $p=1,...,m$, $m\in \bf N$,
$b_p=a_{p+1}$ for each $p=1,...,m-1$, $a_1=0$, $b_m=1$.
Then there exists a rectifiable closed path $\psi $ homotopic to $\eta $
and a neighbourhood $U$ satisfying conditions of Theorem $3.9$ \cite{luoyst}
and such that $\psi \subset U\setminus \sigma (T)$.
For $T\in L_q(X)$ we can define \\
$(i)\quad f(T):=(2\pi )^{-1} (\int_{\psi }f(\zeta )R(\zeta ;T)d\zeta
)M^{-1},$ \\
where convergence is supposed in the weak operator topology.
This integral depends on $f, T$ and it is independent from
$U$, $\psi $, $\eta $, $\gamma $, $M$.
For unbounded QLO $T$ let $A:=-R(a;T)$
and $\Psi : {\bf {\hat H}}\to {\bf {\hat H}}$, $\Psi (z):=(z-a)^{-1}$,
$\Psi (\infty )=0$, $\Psi (a)=\infty $, where $a\in \rho (T)$.
For $f\in {\cal H}_{\infty }(T)$ we define $f(T):=\phi (A)$,
where $\phi \in {\cal H}_{\infty }(A)$ is given by the equation
$\phi (z):=f(\Psi ^{-1}(z))$.
\par {\bf 2.7. Note.} Consider BS $X$ over $\bf H$ as BS $X_{\bf C}$
over $\bf C$, then
\par $(i)$ $X_{\bf C}= X_1\oplus X_2j$,\\
where $X_1$ and $X_2$ are BS over $\bf C$, such that $X_1$ is isomorphic
with $X_2$. The complex conjugation in $\bf C$ induces the complex
conjugation of vectors in $X_m$, where $m=1$ and $m=2$. Each vector
$x\in X$ can be written in the matrix form
\par $(ii)$ $x={{x_1\quad x_2}\choose {-{\bar x}_2\quad {\bar x}_1}}$, \\
where $x_1\in X_1$ and $x_2\in X_2$.
Then each QLO $T$ can be written in the form
\par $(iii)$ $T={{T_1\quad T_2}\choose {-{\bar T}_2\quad {\bar T}_1}}$, \\
where $T_1: X_1\supset {\cal D}(T_1)\to Y_1$,
$T_2: X_1\supset {\cal D}(T_2)\to Y_2$, $T(x)=Tx$
for $s\in \{ q, r \} $, $T(x)=xT$ for $s=l$.
\par $(iv)$ ${\bar T}_mx:=\overline {T_m{\bar x}}$, where $m=1$ or $m=2$.
\par In particular, for the commutator  $[\zeta I,T]$
when $\zeta ={{b\quad 0}\choose {0\quad \bar b}}\in \bf H$,
where $b\in \bf C$, is accomplished the formula
\par $(v)$ $[\zeta I,T]=2(-1)^{1/2} Im(b) {{0\quad T_2}\choose
{-{\bar T}_2 \quad 0}}$, where $Im (b)$ is the imaginary part $b$, \\
$2(-1)^{1/2} Im(b)=(b-{\bar b})$.
\par {\bf 2.8. Theorem.} {\it If $f\in {\cal H}_{\infty }(T)$,
then $f(T)$ does not depend from $a\in \rho (T)$ and    \\
$(i)\quad f(T)=f(\infty )I+(2\pi )^{-1} (\int_{\psi }
f(\lambda ) R(\lambda ;T) d\lambda ) M^{-1}$.}
\par {Proof.} If $a\in \rho (T)$, then $0\ne b=(\lambda
-a)^{-1}$ for $\lambda \ne a$ and $(T-aI)(T-\lambda I)^{-1}=
(bI-A)^{-1}b$, therefore, $I+b^{-1}(T-\lambda I)^{-1}=b(bI-A)^{-1}b-bI$
and $b\in \rho (A)$. If $0\ne b\in \rho (A)$, then $A(bI-A)^{-1}=
(T-\lambda I)^{-1}b^{-1}$ and $\lambda \in \rho (T)$. The point
$b=0\in \sigma (A)$, since $A^{-1}=T-aI$ is unbounded.
Let $a\notin V$, then $U=\psi ^{-1}(V)\supset \sigma (A)$
and $U$ is open in $\bf \hat H$, and $\phi (z):=f(\psi ^{-1}(z))
\in {\cal H}_{\infty }(U)$. In view of Corollary 3.26 \cite{luoyst}
we get $(i)$. 
\par {\bf 2.9. Theorem.} {\it Let $a, b, c, e\in \bf H$,
$f, g \in {\cal H}_{\infty }(T)$. Then
\par $(i)$ $afc+bge\in {\cal H}_{\infty }(T)$ and $afc(T)+bge(T)=
(afc+bge)(T)$;
\par $(ii)$ $fg\in {\cal H}_{\infty }(T)$ and $f(T)g(T)=(fg)(T)$;
\par $(iii)$ if $f$ is decomposed into a convergent series
$f(z)=\sum_k(b_k,z^k)$ in a neighbourhood $\sigma (T)$,
then $f(T)=\sum_k(b_k,T^k)$ on ${\cal D}(T^{\infty })$,
where $b_k=(b_{k,1},...,b_{k,m(k)})$,
$(b_k,z^k):=b_{k,1}z^{k_1}...b_{k,m(k)}z^{k_{m(k)}}$,
$b_{k,j}\in \bf H$;
\par $(iv)$ $f\in {\cal H}_{\infty }(T^*)$ and $f^*(T^*)=(f(T))^*$, where
$f^*(z):=(f(z^*))^*$.}
\par {\bf Proof.} $\bf (i).$ Take $V_f\cap V_g=:V$
and for it we construct $U$, $\eta $ and $\psi $ as in \S 2.6.
Then the first statement follows from $2.6.(i)$.
\par $\bf (ii).$ In view of theorem $3.28$ \cite{luoyst}
the function ${\tilde g}({\tilde z})=:\phi (z)$ belongs to
${\cal H}_{\infty }(T)$ and ${\tilde \phi }({\tilde z})=g(z)$,
where ${\tilde z}=vI-wJ-xK-yL$, $z=vI+wJ+xK+yL$,
$v, w, x, y \in \bf R$, $z\in V_g\subset \bf H$.
Using $\phi (A)$ we consider the case of bounded $T$.
In view of $2.6.(i)$:  $({\tilde y}^*\phi ({\tilde T})
{\tilde h})^{~}=(2\pi )^{-1}y^* M\int_{\psi }d{\tilde \zeta }
R({\tilde \zeta },T)g({\tilde \zeta })h$,
where ${\tilde y}^*{\tilde T}{\tilde h}:=(y^*Th)^{\tilde .}$
and ${\tilde y}^*{\tilde h}:=(y^*h)^{\tilde .}$
for each $y^*\in X^*$ and $h\in X$. Therefore, 
$({\tilde y}^*\phi ({\tilde T}){\tilde h})^{\tilde .}=(2\pi )^{-1}y^*
M\int_{\tilde \psi }d\zeta R(\zeta ,T)g(\zeta )h$, consequently,
$(\phi ({\tilde T}))^{\tilde .}=(2\pi )^{-1}M
\int_{\tilde \psi }d\zeta R(\zeta ,T)g(\zeta )=g(T)$,
since left and right integrals coinside in the space of quaternion
holomorphic functions. The function $fg$ is quaternion holomorphic
on $V$ (see \S \S $2.1$ and $2.12$ \cite{luoyst}).
There exist $\psi _f$ and $\psi _g$ as in \S $2.6$ and contained in
$U\setminus \sigma (T)$, where ${\bar U}\subset V$.
In view of the Fubini theorem there exists
$$(v)\quad f(T)g(T)=(2\pi )^{-2}\int_{\psi _f}\int_{\psi _g}
f(\zeta _1)R(\zeta _1;T)(d\zeta _1)(d\zeta _2)R(\zeta _2,T)g(\zeta _2),$$
where $\zeta _1\in \psi _f$ and $\zeta _2\in \psi _g$.
There are accomplished the identities
$R(\zeta ;T)d\zeta =d_{\zeta } Ln (\zeta I-T)$ and 
$(d\zeta ) R(\zeta ;T)=d_{\zeta} Ln (\zeta I-T)$ for a chosen branch
of $Ln$ (see \S \S $3.7, 3.8$ \cite{luoyst}),
consequently, $R(\zeta _1;T)(d\zeta _1)(d\zeta _2)R(\zeta _2;T)=$
$d_{\zeta _1}d_{\zeta _2} Ln (\zeta _1I-T) Ln (\zeta _2I-T)$
$=(d\zeta _1)R(\zeta _1;T)R(\zeta _2;T)d\zeta _2$.
In view of Lemma $2.2$:
$$(vi)\quad R(a;T)R(b;T)=[R(a;T)-R(b;T)](b-a)^{-1}$$
$$ +R(a;T)[R(b;T),(b-a)I](b-a)^{-1},$$
$$(vii)\quad [R(b;T),(b-a)I]=R(b;T)[T,(b-a)I]R(b;T),$$
since $[(bI-T),(b-a)I]=-[T,(b-a)I]$, where $a, b\in \rho (T)$.
Let in particular $\psi _f$ and $\psi _g$
are contained in the plane ${\bf R}\oplus i\bf R$ in $\bf H$,
where $i, j, k$ are generators of $\bf H$ such that
$i^2=j^2=k^2=-1$, $ij=k$, $jk=i$, $ki=j$.
In view of $2.7.(v)$ and $2.8.(vii)$:
$$(viii)\quad \int_{\psi _f}\int_{\psi _g}f(\zeta _1)R(\zeta _1;T)
[R(\zeta _2;T),(\zeta _2-\zeta _1)I](\zeta _2-\zeta _1)^{-1}d\zeta _1
d\zeta _2 g(\zeta _2)=0,$$
since the branch of $Ln$ can be chosen the same along the axis
$j$ in $\bf H$, in view of the argument principle $3.30$ \cite{luoyst}
it corresponds to the residue $c(\zeta _1-z)^{-1}(\zeta _2-z)^{-1}
b(\zeta _2-\zeta _1)(\zeta _2-z)^{-1}(\zeta _2-\zeta _1)^{-1}$.
Then from $(v,vi,viii)$ it follows:
$$(ix)\quad f(T)g(T)=(2\pi )^{-2}\int_{\psi _f}\int_{\psi _g}
f(\zeta _1)[R(\zeta _1;T)-R(\zeta _2,T)](\zeta _2-\zeta _1)^{-1}
d\zeta _1d\zeta _2 g(\zeta _2).$$
Choose $\psi _g $ such that $|\psi _g(s)-z_0| > |\psi _f(s)-z_0|$
for each $s\in [0,1]$. From the additivity of the integral along the path
and the Fubini theorem:
$$(x) \quad f(T)g(T)=(2\pi )^{-2}\int_{\psi _f}
f(\zeta _1)R(\zeta _1;T)d\zeta _1(\int_{\psi _g}(\zeta _2-\zeta _1)^{-1}
d\zeta _2 g(\zeta _2))$$
$$-(2\pi )^{-2} \int_{\psi _g}
(\int_{\psi _f}f(\zeta _1)d\zeta _1 R(\zeta _2,T) (\zeta _2-\zeta _1)^{-1})
d\zeta _2 g(\zeta _2).$$
In view of Theorems $3.9, 3.24$ \cite{luoyst} the second integral
on the right of $(x)$ is equal to zero, since
$\int_{\psi _f}f(\zeta _1)d\zeta _1 R(\zeta _2;T)
(\zeta _2-\zeta _1)^{-1}=0$, the first integral produces
$$f(T)g(T)=(2\pi )^{-1} (\int_{\psi _f}
f(\zeta _1) R(\zeta _1;T) d\zeta _1 g(\zeta _1))M^{-1}$$
$$=(2\pi )^{-1} (\int_{\psi _f}
f(\zeta )g(\zeta ) d Ln (\zeta I-T)) M^{-1},$$
where $\zeta \in \psi _f$.
\par $\bf (iii)$ follows from the application of $\bf (i,ii)$
by induction and convergence of the series in the strong operator
topology.
\par $\bf (iv).$ In view of Lemma $2.5$
$\sigma (T^*)=(\sigma (T))^{\tilde .}$,
then $f\in {\cal H}_{\infty }(T^*)$.
Since $(f(T))^*y^*=y^*\circ f(T)$ for each $y^*\in X^*$, then due to
Lemma 2.5 $R(\zeta ^*;T^*)=(R(\zeta ;T))^*$, consequently,
$(f(T))^*y^*=(2\pi )^{-1}(M^{-1})^*\int_{\psi }
(d\zeta ^*)R(\zeta ^*;T^*) (f(\zeta ))^*y^*$, where
$(f(\zeta ))^*y^*:=y^*\circ f(\zeta )$. If $f(\zeta )$
is represented by the series converging in the ball:
$f(\zeta )=\sum_n (a_n,\zeta ^n)$, then
$f(\zeta ^*)=\sum_na_{n,1}\zeta ^{*n_1}...a_{n,m(n)}\zeta ^{*n_{m(n)}}$,
consequently,
$[f(\zeta ^*)]^*=\sum_n\zeta ^{n_{m(n)}}a_{n,m(n)}^*...\zeta ^{n_1}
a_{n,1}^*$ and \\
$(f(T))^*=(2\pi )^{-1}(M^{-1})^*\int_{\psi }(d\zeta ^*)R(\zeta ^*;T^*)
f^*(\zeta ^*)=f^*(T^*)$.
\par {\bf 2.10. Theorem.} {\it Let QLO $T$ be bounded,
$f\in {\cal H}(T)$, then
$f(\sigma (T))=\sigma (f(T)).$}
\par {\bf 2.11. Theorem.} {\it Let $f\in {\cal H}_{\infty }(T)$,
let also $f(U)$ be open for some open $U\subset dom (f)\subset \bf \hat H$,
$g\in {\cal H}_{\infty }(T)$ and $f(U)\supset \sigma (T)$,
$dom (g)\supset f(U)$, then $F:=g\circ f\in {\cal H}(T)$ and
$F(T)=g(f(T)$.}
\par {\bf Proof} follows from \S 2.9 analogously to the case
of the field $\bf C$
with the help of Theorems $2.16, 3.10$, Corollary $2.13$ \cite{luoyst},
since $F\in {\cal H}_{\infty }(T)$.
\par {\bf 2.12. Definition and Note.} Let $\cal A$ be BS and an algebra
over $\bf H$ with the unity $e$ having the properties:
$|e|=1$ and $|xy|\le |x| |y|$ for each $x$ and $y\in \cal A$,
then $\cal A$ is called the Banach algebra (BA) or $C$-algebra
(over $\bf H$). BA $\cal A$ is called quasi-commutative (QC), if
there exists a commutative algebra ${\cal A}_{0,0}$ over $\bf R$
such that ${\cal A}$ is isomorphic with the algebra
$\{ T: T={{A\quad B}\choose {{-\bar B}\quad {\bar A}}}:$
$A, B\in {\cal A}_0 \} $, where ${\cal A}_0:=\{ A:$
$A=A_0+A_1\bf i$; $A_0\in {\cal A}_{0,0},$ $A_1\in {\cal A}_{0,0} \} $,
${\bar A}:=A_0-A_1\bf i$, ${\bf i}:=(-1)^{1/2}$.
\par Consider $X$ over $\bf R$: $\quad X_{\bf R}=X_ee\oplus
X_ii\oplus X_jj\oplus X_kk$, where $X_e$, $X_i$, $X_j$
and $X_k$ are pairwise isomorphic BS over $\bf R$.
Then ${\cal A}={\cal A}_0e\oplus {\cal A}_ii\oplus {\cal A}_jj\oplus
{\cal A}_kk$, where ${\cal A}_0$, ${\cal A}_i$, ${\cal A}_j$
and ${\cal A}_k$ are algebras over $\bf R$. Multiplying
$\cal A$ on $S\in \{ e,i,j,k \} $, we get automorphisms of $\cal A$,
consequently, ${\cal A}_0$, ${\cal A}_i$, ${\cal A}_j$ and ${\cal A}_k$
are pairwise isomorphic.
\par {\bf 2.13. Definitions and Notes.}
BA $\cal A$ is supplied with the involution, when there exists an operation
$*: {\cal A}\ni T\mapsto T^*\in \cal A$ such that
$(T^*)^*=T$, $(T+V)^*=T^*+V^*$,
$(TV)^*=V^*T^*$, $(\alpha T)^*=T^*\tilde \alpha $
for each $\alpha \in \bf H$.
\par If BA $\cal A$ (over $\bf H$) has a subalgebra ${\cal A}_{0,0}$
(over $\bf R$), then
$T^*={{A^*\quad -{\bar B}^*}\choose {B^*\quad {\bar A}^*}}.$
\par An element $x\in \cal A$ is called regular, if there exists
$x^{-1}\in \cal A$. In the contrary case it is called singular.
Then the spectrum $\sigma (x)$ for $x$ is defined as the set of all
$z\in \bf H$, for which $ze-x$ is singular, his spectral radius
is $|\sigma (x)|:=\sup_{z\in \sigma (x)} |z|$.
The resolvent set is defined as $\rho (x):= \{ z\in {\bf H}:$
$ze-x$ $\mbox{is regular} \} $ and the resolvent is
$R(z;x):=(ze-x)^{-1}$ for each $z\in \rho (x)$.
\par {\bf 2.14. Lemma.} {\it A spectrum $\sigma (x)$
of an element $x\in \cal A$ is a non-void compact subset
in $\bf H$. Its resolvent $x(z):=R(z;x)$ is quaternion holomorphic
on $\rho (x)$, $x(z)$ converges to zero when $|z|\to \infty $ and
\par $x(z)-x(y)=x(z)(y-z)x(y)$ for each $y, z\in \rho (x)$.}
\par {\bf Proof} follows from $(ze-x) x(z)x(y)=x(y)$,
$x(z)x(y) (ye-x)=x(z)$, $(ze-x)(x(z)-x(y))(ye-x)=(ye-x)-(ze-x)$
$=(y-z)e$, consequently, $x(z)-x(y)=R(z;x)(y-z)R(y;x)$
$=x(z)(y-z)x(y)$. Therefore, $x(z)$ is continuous by $z$ on
$\rho (x)$ and there exists $(\partial [x(z+y)
x^{-1}(y)]/\partial z).h=-x(y)h$ for each $h\in \bf H$.
For each marked point $y$ the term $x^{-1}(y)$ is constant on
$\cal A$, moreover, $(\partial [x(z+y)x^{-1}(y)]/\partial {\tilde z})=0$,
consequently, $x(z) \in {\cal H} (\rho (x))$.
The second statement follows from the consideration of the complexification
${\bf C}\otimes \cal A$.
\par {\bf 2.15. Theorem.} {\it Let $\cal B$ be a closed ideal
in QCBA $\cal A$. The quotient algebra
${\cal A}/\cal B$ is isometrically isomorphic with $\bf H$
if and only if $\cal B$ is maximal.}
\par The {\bf proof} is analogous to the case of algebras over $\bf C$
due to the definition of QCBA.
\par {\bf 2.16. Definitions.} A $C^*$-algebra $\cal A$ over $\bf H$
is a BA over $\bf H$ with the involution $*$ such that
$|x^*x|=|x|^2$ for each $x\in \cal A$.
\par A scalar product on a linear space $X$ over $\bf H$
(that is, linear relative to the right and left multiplications
separately on scalars from $\bf H$) is the biadditive $\bf R$-bilinear
mapping $<*;*>: X^2\to \bf H$ such that
\par $(1)$ $<x;x>=\alpha _0e$, where $\alpha _0\in \bf R$;
\par $(2)$ $<x;x>=0$ if and only if $x=0$;
\par $(3)$ $<x;y>=<y;x>^{\tilde .}$ for each $x, y \in X$;
\par $(4)$ $<x+z;y>=<x;y>+<z;y>$;
\par $(5)$ $<xa;yb>={\tilde a}<x;y>b$ for each
$x, y, z\in X$, $a, b\in \bf H$.
\par If $X$ is complete relative to the norm topology
\par $(6)$ $|x|:=<x;x>^{1/2}$, then $X$ is called the quaternion
Hilbert space (HS).
\par {\bf 2.17. Lemma.} {\it BA $L_q(X)$ on HS $X$ with an involution:
\par $(1)$ $<Tx;y>=:<x;T^*y>$ for each $x, y \in X$ \\
is a $C^*$-algebra.}
\par {\bf 2.18. Lemma.} {\it If $\cal A$ is a QC 
$C^*$-algebra, then $|x^2|=|x|^2$, $|x|=|x^*|$ and $I^*=I$,
where $I$ is the unity in $\cal A$.}
\par {\bf Proof.} Each vector $x\in \cal A$
we represent in the form: $x=x_ee+x_ii+x_jj+x_kk$.
Then $x^*=x_e^*e-x_i^*i-x_j^*j-x_k^*k$, since
$(x_mS_m)^*=(-1)^{\kappa (S_m)}x_mS_m$, where
$S_m\in \{ e, i, j, k \}$ for each $m\in \{ e,i,j,k \}$,
$\kappa (e)=0$, $\kappa (i)=\kappa (j)=\kappa (k)=1$.
Therefore, $\quad [x,x^*]=0$ and $|x_m^2|=|x_m|^2$.
Then $|x|^2=|x_e|^2+|x_i|^2+|x_j|^2+|x_k|^2$ and
$|x^2|^2=|(x^2)^*x^2|=|(x^*)^2x^2|=|(xx^*)(xx^*)|=|x|^4$,
consequently, $|x^2|=|x|^2$. Since $I=I_e$, then $I^*=I_e^*=I_e=I$.
\par {\bf 2.19. Definition.} A homomorphism $h: {\cal A}\to
\cal B$  of $C^*$-algebras $\cal A$ and $\cal B$ preserving involutions:
$h(x^*)=(h(x))^*$ is called a $*$-homomorphism.
If $h$ is a bijective $*$-homomorphism of $\cal A$ on $\cal B$,
then $h$ is called a $*$-isomorphism, $\cal A$ and $\cal B$ are called
$*$-isomorphic. By $\sigma ({\cal A})$ is denoted a structural space
for $\cal A$ and it is aslo called a spectrum of $\cal A$.
A structural space is defined analogously to the complex case
with the help of Theorem 2.15.
\par {\bf 2.20. Proposition.} {\it For HS $X$ spaces
$L_l(X,{\bf H})$ and $L_r(X,{\bf H})$ are isomorphic with $X$, for BS $X$
$L_q(X)$ is isomorphic with $L_l(X^2)$, $L_q({\bf H})=\bf H^4$,
there exists a bijection between a family of all QLO $T$ on
${\cal D}(T)\subset X$ a family of all LLO $V$ on ${\cal D}(V)\subset X^2$.}
\par {\bf Proof.} Let $L_q(X)\ni \alpha =\alpha _ee+
\alpha _ii+\alpha _jj+\alpha _kk$.
Since $S\alpha S\in L_q(X)$ for each
$S\in \bf H$, there exist quaternion constants $S_{m,l,n}$ such that
\par $(1)\quad \alpha _m(x)m=
\sum_nS_{m,1,n}\alpha (x)S_{m,2,n}$ for each $m\in \{ e,i,j,k \} $, \\
where $S_{m,1,n}=\gamma _{m,n}S_{m,2,n}$ with
$\gamma _{m,n}=(-1)^{\phi (m,n)}/4\in \bf R$,
$\phi (m,n)\in \{ 1, 2 \} $, $S_{m,l,n}\in {\bf R}n$
for each $n\in \{ e,i,j,k \} $ (see \S \S 3.7, 3.28 \cite{luoyst}).
Applying for $x$ the decomposition from \S 2.12,
due to $(1)$ we get $4\times 4$-block form
of operators over $\bf R$ and the isomorphism of $L_q(X)$ with $L_l(X^2)$.
\par {\bf 2.21. Note.} For LLO (over $\bf H$) the notions of point
$\sigma _p(T)$, continuous $\sigma _c(T)$ and residual $\sigma _r(T)$
spectra are defined analogously to the case over the field $\bf C$,
due to 2.20 these notions spread on QLO.
\par {\bf 2.22. Theorem.} {\it QC $C^*$-algebra
is isometrically $*$-isomorphic with the algebra $C(\Lambda ,{\bf H})$
of all continuous $\bf H$-valued functions on its spectrum $\Lambda $.}
\par The {\bf proof} follows from the fact, that the mapping
$x\mapsto x(.)$ from $\cal A$ into $C(\Lambda ,{\bf H})$
is the $*$-homomorphism, where $x({\cal M})$ is defined by the
equality $x+{\cal M}=x({\cal M})+{\cal M}$ for each maximal
ideal $\cal M$. Let $x(\lambda )=\alpha _ee+
\alpha _ii+\alpha _jj+\alpha _kk$, then $x^*(\lambda )=
\beta _ee+\beta _ii+\beta _jj+\beta _kk$, where $\alpha _e,...,
\beta _k\in \bf R$.
There exists the decomposition for $X:=C(\Lambda ,{\bf H})$ from \S 2.12
with $X_e=C(\Lambda ,{\bf H})$
and $x_m=\sum_nS_{m,1,n}zS_{m,2,n}{\tilde m}$ for each $z\in \bf H$,
where $z=x_ee+x_ii+x_jj+x_kk$, $x_m\in \bf R$, $m, n\in \{ e,i,j,k \} $
(see Proposition 2.20).
Therefore, it can be applied the Stone-Weiestrass theorem for $\bf H$-valued
functions. If $\lambda _1\ne \lambda _2$ are two maximal ideals in
$\Lambda $, then $y(\lambda _1)\ne y(\lambda _2)$ for each
$y\in \lambda _1\setminus
\lambda _2$. Consequently, the algebra of functions $x(.)$ coinsides
with $C(\Lambda ,{\bf H})$.
\par {\bf 2.23. Definition.} Let $X$ be HS over $\bf H$ and $\cal B$
be a $\sigma $-algebra of Borel subsets of a Hausdorff topological space
$\Lambda $. Consider a mapping $\hat E$ defined on ${\cal B}\times X^2$
and defining a unique $X$-projection-valued spectral measure
$\hat E$ such that
\par $(i)$ $<{\hat E}(\delta )x;y>={\hat \mu }(\delta ;x,y)$
is a regular (non-commutative) $\bf H$-valued measure for each $x, y \in X$,
where $\delta \in \cal B$. By our definition this means, that
\par $(1)\quad {\hat \mu }(\delta ;x,y)=
{\hat \mu }(x,y).\chi _{\delta }$ and
\par $(2)\quad {\hat \mu }(x,y).f:=
\sum_{m,l,n}\int_{\Lambda }S_{m,1,n}f(\lambda )S_{m,2,n}{\tilde m}l
\mu _{m,l}(d\lambda ;x,y)$, \\
where $\chi _{\delta }$ is the characteristic function for
$\delta \in \cal B$, $S_{m,p,n}\in {\bf R}n$ are the same as in \S 2.19,
$\mu _{m,l}$ is a regular real-valued measure,
$m,n,l \in \{ e,i,j,k \} $, $p=1$ or $p=2$,
$f$ is an arbitrary $\bf H$-valued function on $\Lambda $,
which is $\mu _{m,l}$-integrable for each $m, l$;
\par $(3)\quad ({\hat E}_S(\delta ))^*=(-1)^{\kappa (S)}
{\hat E}_S(\delta ),$ where $({\hat E}_S(\delta ).e)x:=
({\hat E}_S(\delta ))x=({\hat E}(\delta ).S)x$,
$S=cm$, $m\in \{ e,i,j,k \} $, $c=const\in \bf R$, $x\in X$;
\par $(4)\quad {\hat E}_{bS}=b{\hat E}_S$ for each $b\in \bf R$
and each pure vector $S=cm$;
\par $(5)\quad {\hat E}_{S_1S_2}(\delta \cap \gamma )=
{\hat E}_{S_1}(\delta ){\hat E}_{S_2}(\gamma )$ for each
pure quaternion vectors $S_1$ and $S_2$ and each $\delta , \gamma
\in \cal B$. Though from $(3,4)$ it follows, that
${\hat E}(\delta ).\lambda =
\lambda _e{\hat E}_e(\delta )+\lambda _i{\hat E}_i(\delta )
+\lambda _j{\hat E}_j(\delta )+\lambda _k{\hat E}_k(\delta ) =:
{\hat E}_{\lambda }(\delta )$, but in general it may be
$({\hat E}(\delta ).\lambda )x\ne ({\hat E}(\delta ))\lambda x$,
where $\lambda _e,...,\lambda _k\in \bf R$.
\par {\bf 2.24. Theorem.} {\it Each QC $C^*$-algebra
$\cal A$ contained in $L_q(X)$ for HS $X$ over $\bf H$
is isometrically $*$-equivalent with the algebra $C(\Lambda ,{\bf H})$,
where $\Lambda $ is its spectrum. Moreover, each isometrical $*$-isomorphism
$f\mapsto T(f)$ between $C(\Lambda ,{\bf H})$ and $\cal A$
defines a unique $X$-projection-valued spectral measure $\hat E$ on
${\cal B}(\Lambda )$ such that
\par $(i)$ $<{\hat E}(\delta )x;y>={\hat \mu }(\delta ;x,y)$
is the regular $\bf H$-valued measure for each $x, y \in X$,
where $\delta \in \cal B$;
\par $(ii)$ ${\hat E}_{S_1}(\delta ).T(S_2f)=(-1)^{\kappa (S_1)+
\kappa (S_2)}T(S_2{\hat E}_{S_1}(\delta ).f)$
for each $f\in C(\Lambda ,{\bf R})$, $\delta \in \cal B$
and pure quaternion vectors $S_1$, $S_2$;
\par $(iii)\quad T(f)=\int_{\Lambda }{\hat E}(d\lambda ).f(\lambda )$
for each $f\in C(\Lambda ,{\bf H})$,
moreover, $\hat E$ is $\sigma $-additive in the strong operator
topology.}
\par {\bf Proof.} We mention, that $\Lambda $ is compact.
There exists a decomposition $C(\Lambda ,{\bf H})$ as in \S 2.21.
Each $\psi \in C^*_q(\Lambda ,{\bf H})$
has a decomposition $\psi (f)=\psi _e(f)e+\psi _i(f)i+\psi _j(e)j
+\psi _k(f)k$, where $f\in C(\Lambda ,{\bf H})$ (see $2.20.(1)$).
Moreover, $\psi _l(f)=\psi _l(f_ee)+\psi _l(f_ii)+\psi _l(f_jj)+
\psi _l(f_kk)$, where $f_m\in C(\Lambda ,{\bf R})$,
$m, l \in \{ e,i,j,k \} $. Then
\par $(1)\quad \psi (f)=\sum_{m,n,l}\psi _l(S_{m,1,n}fS_{m,2,n})l$, \\
where $m, n, l \in \{ e,i,j,k \} $. In view of the Riesz representation
theorem $IV.6.3$ \cite{danschw}: $\psi _l(gm)=\int_{\Lambda }g(\lambda )
\mu _{m,l}(d\lambda )$ for each $g\in C(\Lambda ,{\bf R})$,
where $\mu _{m,l}$ is a $\sigma $-additive real-valued measure.
The accomplishment of the componentwise integration of matrix-valued
functions gives
\par $(2)\quad \psi (f)=\sum_{m,n,l} \int_{\Lambda }
S_{m,1,n}f(\lambda )S_{m,2,n}{\tilde m}l\mu _{m,l}(d\lambda )$.
For $\psi (f):=<T(f)x;y>$ for each $f\in C(\Lambda ,{\bf H})$
and marked $x, y\in X$ from $(2)$ it follows, that
\par $(3)\quad <T(f)x;y>=\sum_{m,n,l} \int_{\Lambda }
S_{m,1,n}f(\lambda )S_{m,2,n}{\tilde m}l\mu _{m,l}(d\lambda ;x,y)$,
since $|<T(f)x;y>| \le |f| |x| |y|$, consequently,
$\mu _{m,l}(\delta ;xa,yb)=a\mu _{m,l}(\delta ;x,y)b$
for each $a, b\in \bf R$, moreover,
\par $(4)\quad \sup_{\delta \in \cal B} (\sum_l
|\sum_m z_mm\mu _{m,l}(\delta ;x,y)|^2)^{1/2}\le |z| |x| |y|$
for each $z=z_ee+z_ii+z_jj+z_kk\in \bf H$, since $|l|=1$.
\par From $(3)$ it follows, that $\mu _{m,l}(\delta ;x,y)$
is $\bf R$-bihomogeneous and biadditive by $x, y$.
if $f(\lambda )\in {\bf R}m$ for $\mu $-a.e. $\lambda \in \Lambda $
for some $m\in \{ e,i,j,k \} $, then
$T(f)=T((-1)^{\kappa (m)}{\tilde f})=(-1)^{\kappa (m)}T(f)^*$, consequently,
$<T(f)x;y>=(-1)^{\kappa (m)}<T(f)y;x>^{\tilde .}$.
Then $\mu _{m,l}(\delta ;x,y)=(-1)^{\kappa (m)+\kappa (l)}\mu _{m,l}
(\delta ;y,x)$ for each $m, l\in \{ e,i,j,k \} $, $x, y\in X$.
\par {\bf 2.25. Definition.} An operator $T$ on a quaternion HS
$X$ is called normal, if $TT^*=T^*T$; $T$ unitary, if
$TT^*=I$ and $T^*T=I$; $T$ symmetrical, if $<Tx;y>=<x;Ty>$ for each
$x, y \in {\cal D}(T)$, $T$ self-adjoint, if $T^*=T$.
Henceforth, for $T^*$ it is supposed, that ${\cal D}(T)$ is dense
in $X$.
\par {\bf 2.26. Lemma.} {\it An operator $T\in L_q(X)$
is normal if and only if a minimal ($\bf H$-)subalgebra $\cal A$
in $L_q(X)$ containing $T$ and $T^*$ is QC.}
\par {\bf Proof.} Let $T$ be normal, then on
$X=X_ee\oplus X_ii\oplus X_jj\oplus X_kk$ it can be represented in the form
$T=T_eE+T_ii+T_jj+T_kk$, where $Range (T_m)\subset X_m$
and $T_m\in \cal A$ for each $m\in \{ e,i,j,k \} $. In the block form
$T={{A\quad B}\choose {-{\bar B}\quad {\bar A}}}$ and $x={{x_1\quad
x_2}\choose {-{\bar x}_2\quad {\bar x}_1}}$, where $x_1=x_ee+x_ii\in
X_ee\oplus X_ii$, $x_2j=x_jj+x_kk\in X_jj\oplus X_kk$,
${\bar x}_1=x_ee-x_ii$, ${\bar A}x_1:=A{\bar x}_1$.
Then $T^*(x_mm)=(-1)^{\kappa (S_m)}T(x_mm)$ for each
$m\in \{ e,i,j,k \} $ and
$T^*={{A^*\quad -{\bar B}^*}\choose {B^*\quad {\bar A}^*}}$.
Therefore, $<mTmx;mTmy>=<mT^*mx;mT^*my>$ for each $m\in \bf H$
with $|m|=1$, consequently, $(mTm)(mTm)^*=(mTm)^*(mTm)$.
The space $X$ is isomorphic with $l_2(\upsilon ,{\bf H})$, in which
$<x;y>=\sum_{b\in \upsilon }\mbox{ }^b{\tilde x}\mbox{ }^by$,
where $\upsilon $ is a set, $x= \{ \mbox{ }^lx:
\mbox{ }^lx\in {\bf H}, l\in \upsilon \} \in l_2(\upsilon ,{\bf H})$.
Then in $L_q(l_2(\upsilon ,{\bf H}))$ is accomplished $T^*=\tilde T$,
moreover, ${\bar A}^*=A$ and ${\bar B}^*=B$.
Therefore, $TT^*=T^*T$ gives $AA^*=A^*A$, also an automorphism
$j: X\to X$ and the equality $(mTm)(mTm)^*=(mTm)^*(mTm)$
with $m=\theta $, $m=(i+j)/2$, $m=(i+k)/2$, $m=(i+j)\theta /2$,
$m=(i+k)\theta /2$, $\theta =\exp (\pi {\bf i}/4){{1\quad 0}
\choose {0\quad {-1}}}$ leads to the pairwise commuting
$\{ T_e,T_i,T_j,T_k \} $.
\par Vice versa, if $\cal A$ is quasi-commutative, then
$\{ T_m: m=e,i,j,k \} $ are pairwise commuting, consequently,
$TT^*=T^*T$.
\par {\bf 2.27. Lemma.} {\it Let $T$ be a symmetrical operator and
$a\in {\bf H}\setminus {\bf R}e$, then there exists
$R(a;T)$ and $|x|\le 2 |R(a;T)x| / |a-{\tilde a}|$ for each
$x\in {\cal D}(T)$. Let $T$ be a closed operator, then the sets
$\rho (T)$, $\sigma _p(T)$, $\sigma _c(T)$ and $\sigma _r(T)$
do not intersect and their union is the entire $\bf H$. For a self-adjoint
QLO $T$ $\sigma (T)\subset {\bf R}e$, moreover, $R(a;T)^*=R(a^*;T)$.}
\par {\bf 2.28. Theorem.} {\it For a self-adjoint QLO $T$
there exists a uniquely defined regular countably-additive
self-adjoint spectral measure $\hat E$ on ${\cal B}({\bf H})$,
${\hat E}|_{\rho (T)}=0$ such that \\
$(a)$ ${\cal D}(T):= \{ x: x\in X; \int_{\sigma (T)}<({\hat E}(dz).z^2)x;x>
<\infty \} $ and \\
$(b)$ $Tx=\lim_{n\to \infty }\int^n_{-n}({\hat E}(dz).z)x$,
$x\in {\cal D}(T)$.}
\par {\bf Proof.} In view of Proposition 2.20 and Equality 2.16.(5)
the space ${\cal D}(T)$ is $\bf H$-linear.
Let us use Lemma 2.27, the proof of which is analogous to the case
over the field $\bf C$, also take a marked element
$q\in \{ i,j,k \} $, then there is
$h(z):=(q-z)^{-1}$ the homeomorphism of the sphere $S^3:=\{ z\in {\bf H}:
|z|=1 \} $ and for $A:=(q-z)R(z;T)(q-z)+(q-z)I$ for each
$z\in \rho (T)\setminus \{ q \} $ is accomplished the equality
$(hI-R(q;T))A=I$. If $z=q$, then $h=\infty $,
consequently, $h\notin \sigma (R(q;T))$.
Let $0\ne h\in \rho (R(q;T))$, then there exists $B:=
R(q;T)A$, where $A:=(hI-R(z;T))^{-1}$, consequently, $B$ is bijective,
${\cal R} (B)={\cal D}(T)$ and $(zI-T)B=(z-q)I$, that is,
$z\in \rho (T)$. For $h=0\in \rho (R(q;T))$ the operator
$R(q;T)^{-1}=(hI-T)$ is the bounded everywhere defined operator
and this case is considered in Theorem 2.24. For each
$\delta \in {\cal B}({\bf H})$ we
put ${\hat E}(\delta ):={\hat E}^1(h(\delta ))$, where
${\hat E}^1$ is the decomposition of the identity for the normal operator
$R(q;T)$, then the end of the proof is analogous to that of Theorem
XII.2.3 \cite{danschw}.
\par {\bf 2.29. Note and Definition.} A unique spectral measure,
related with a self-adjoint QLO $T$ is called the decomposition
of the identity for $T$. For $\bf H$-valued Borel function
$f$ defined $\hat E$-almost everywhere $f(T)$ is defined
by the relations:
${\cal D}(f(T)):=\{ x:$ $\mbox{there exists}$ $\lim_nf_n(T)x \} $,
where $f_n(z):=f(z)$ for $|f(z)|\le n$; $f_n(z):=0$ for $|f(z)|>n$;
$f(T)x:=\lim_nf_n(T)x$, $x\in {\cal D}(f(T))$, $n\in \bf N$.
\par {\bf 2.30. Theorem.} {\it Let $\hat E$ be a decomposition
of the identity for a self-adjoint QLO $T$ and $f$ from \S 2.28.
Then $f(T)$ is a closed QLO defined on an everywhere dense domain,
moreover: \\
$(a)$ ${\cal D}(f(T))=\{ x: \int^{\infty }_{-\infty }|f(z)|^2
<{\hat E}(dz)x;x> <\infty \} $; \\
$(b)$ $<f(T)x;y>=\int^{\infty }_{-\infty }<{\hat E}(dz).f(z)x;y>$,
$x\in {\cal D}(f(T))$; \\
$(c)$ $|f(T)x|^2=\int^{\infty }_{-\infty }|f(z)|^2<{\hat E}(dz)x;x>$,
$x\in {\cal D}(f(T))$; \\
$(d)$ $f(T)^*={\tilde f}(T)$; $(e)$ $R(q;T)=\int^{\infty }_{-\infty }
{\hat E}(dz).(q-z)$, $q\in \rho (T)$.}
\par {\bf Proof.} Take $f_n$ from \S 2.29
and $\delta _n:= \{ z: |f(z)|\le n \} $. Then $|f(T)x|^2=
\lim_n|f_n(T)x|^2=\int^{\infty}_{-\infty }|f(z)|^2
<{\hat E}(dz)x;x>$ for each $x\in {\cal D}(f(T))$, from this it follows
$(c)$, a closedness of $f(T)$ and $(a)$ can be verified analogously
to the complex case.
A non-commutative measure $\hat \mu $ on the algebra $\Upsilon $
of subsets of the set $\cal S$ corresponds to QLO with values
in $\bf H$ and due to Proposition 2.20 it is characterized completely by 
$\bf R$-valued measures $\mu _{m,n}$ such that $\mu _{m,n}(f_m)=
{\hat \mu }(f_m)\tilde n$ for each $\hat \mu $-integrable
$\bf H$-valued function $f$ with components $f_m$, where $m,n
\in \{ e,i,j,k \} $. Then it can be defined the variation
$v({\hat \mu },U):=\sup_{W_l\subset U} \sum_l|{\hat \mu }
(\chi _{W_l})|$ by all $\{ W_l \} $ finite disjunctive subsets
$W_l\in \Upsilon $ in $U$ with $\bigcup_lW_l=U$.
If $\hat \mu $ is bounded, then it is QLO with bounded variation
$v({\hat \mu },{\cal S})\le 16 \sup_{U\in \Upsilon }
|{\hat \mu }(\chi _U)|$, moreover, $v({\hat \mu },*)$ is additive on
$\Upsilon $. A function $f$ we call $\hat \mu $-measurable,
if each $f_m$ is $\mu _{m,n}$-measurable for each $n$ and $m
\in \{ e,i,j,k \} .$ The space of all ${\hat \mu }$-measurable
$\bf H$-valued functions $f$ with
$v({\hat \mu },|f|^p)^{1/p}=:|f|_p<\infty $
we denote by $L^p({\hat \mu })$ for $0<p<\infty $,
$L^{\infty }({\hat \mu })$
denotes the space of all $f$ for which there exist
$|f|_{\infty }:=ess_{v({\hat \mu },*)}-\sup |f|<\infty .$
In details we write $L^p({\cal S},{\Upsilon },{\hat \mu },{\bf H})$
instead of $L^p({\hat \mu })$.
A subset $V$ in $\cal S$ we call $\cal \mu $-zero-set,
if $v^*({\hat \mu },V)=0$, where $v^*$ is an extension
of the total variation $v$ by the formula
$v^*({\hat \mu },A):=\inf_{{\Upsilon }\ni F\supset
A}v({\hat \mu },F)$ for $A\subset \cal S$. A non-commutative measure
$\hat \lambda $ on $\cal S$ we call absolutely continuous relative to
$\hat \mu $, if $v^*({\hat \lambda },A)=0$ for each subset
$A\subset \cal S$ with $v^*({\hat \mu },A)=0$.
A measure $\hat \mu $ we call positive, if each $\mu _{m,n}$
is non-negative and $\sum_{m,n}\mu _{m,n}$ is positive. The usage
of components $\mu _{m,n}$ and the classical Radon-Nikodym theorem (
see Theorems III.10.2,10.7) lead to the following its non-commutative
variants.
\par $(i).$ If $({\cal S},{\Upsilon },{\hat \mu })$
is a space with a $\sigma $-finite positive non-commutative
$\bf H$-valued measure $\hat \mu $, $\hat \lambda $ is absolutely
continuous relative to $\hat \mu $ and it is a finite
non-commutative measure defined on $\Upsilon $, then there exists
a unique $f\in L^p({\cal S},{\Upsilon },{\hat \mu },{\bf H})$
such that ${\hat \lambda }(U)={\hat \mu }(f\chi _U)$ for each
$U\in \Upsilon $, moreover, $v({\hat \mu },{\cal S})=|f|_1$.
\par $(ii).$ If $({\cal S},{\Upsilon },{\hat \mu })$
is a space with a finite non-commutative $\bf H$-valued measure $\hat \mu $,
$\hat \lambda $ is absolutely continuous relative to $\hat \mu $
and it is a non-commutative measure defined on $\Upsilon $, then there
exists a unique $f\in L^1({\hat \mu })$ such that ${\hat \lambda }
(U)={\hat \mu }(f\chi _U)$ for each $U\in \Upsilon $.
In view of $(ii)$ there exists a Borel measurable function $\phi $
such that ${\hat \nu }(\delta ):={\hat \mu }_{x,y}(\phi \chi_{\delta })=
<{\hat E}(\phi \chi _{\delta })x;y>$ for each $\delta \in {\cal
B}({\bf R})$. In view of $(i)$ $|\phi (z)|=1$ $\hat \nu $-almost every.
Consider $f^1(z):=|f(z)|\phi (z)$, then in view of $(a)$
${\cal D}(f^1(T))={\cal D}(f(T))$ and $<f^1(T)x;y>=
\int^{\infty }_{-\infty }|f(z)|{\hat \nu }(dz)$.
Therefore, $<f(T)x;y>=\lim_n\int_{\delta _n}<{\hat E}(dz).f(z)x;y>=$
$\int^{\infty }_{-\infty }<{\hat E}(dz).f(z)x;y>$ and from this
it follows $(b)$.
\par $(d)$. From ${\hat E}_S^*=(-1)^{\kappa (S)}{\hat E}_S$
for each $S=cs$, $0\ne c\in \bf R$, $s\in \{ e,i,j,k \} $,
it follows, that ${\hat E}.{\tilde f}={\hat E}^*.f$.
Take $x,y\in {\cal D}({\tilde f}(T))={\cal D}(f(T))$, then
$<{\tilde f}(T)x;y>=\int^{\infty }_{-\infty }
<{\hat E}(dz).{\tilde f}(z)x;y>=<x;f(T)y>$, consequently,
${\tilde f}(T)\subset f(T)^*$. If $y\in {\cal D}(f(T)^*)$, then
for each $x\in X$ and $m\in \bf N$: ${\tilde f}_m(T)y={\hat E}(\delta _m)
.f(T)^*y$ converges to $f(T)^*y$ while $m\to \infty $, consequently,
$y\in {\cal D}({\tilde f}(T))$.
In view of Theorem 2.24 Statement $(e)$ follows from the fact, that
$\mbox{ }_n{\hat E}(\delta ):={\hat E}(\delta _n\cap \delta )$
it is the decomposition of the identity for the restriction
$T|_{X_n}$, where $X_n:={\hat E}(\delta _n)X$.
\par {\bf 2.31. Theorem.} {\it A bounded normal operator
$T$ on a quaternion HS is unitary, Hermitian or positive definite
if and only if $\sigma (T)$ is contained in $S^3:=
\{ z\in {\bf H}: |z|=1 \} $, $\bf R$ or in $[0,\infty )$ respectively.}
\par {\bf Proof.} In view of Theorem 2.24 the equality $T^*T=TT^*=I$
is equivalent to $z{\tilde z}=1$ for each $z\in \sigma (T)$.
If $\sigma (T)\subset [0,\infty )$, then $<Tx;x>=\int_{\sigma (T)}
<{\hat E}(dz).zx;x>\ge 0$ for each $x\in X$. The final part of the proof
is analogous to the complex case, using the technique given above.
\par {\bf 2.32. Definition.} The family $ \{ T(t): 0\le t\in {\bf R} \} $
of bounde QLO in $X$ is called a strongly continuous semigroup,
if $(i)$ $T(t+q)=T(t)T(q)$ for each $t, q\ge 0$;
$(ii)$ $T(0)=I$; $(iii)$ $T(t)x$ is the continuous function by
$t\in [0,\infty )$ for each $x\in X$.
\par {\bf 2.33. Theorem.} {\it For each strongly continuous
semigroup $ \{ U(t): 0\le t\in {\bf R} \} $ of unitary QLO
in HS $X$ over $\bf H$ there exists a unique self-adjoint QLO
$B$ in $X$ such that $U(t)=\exp (t{\bf i}B)$, where ${\bf i}=(-1)^{1/2}$.}
\par {\bf Proof.} If $\{ T(t): 0\le t \} $ is a semigroup continuous
in the uniform topolgy, then due to Theorem
VIII.1.2 \cite{danschw} and Proposition 2.20 there exists a bounded
operator $A$ in $X$ such that $T(t)=\exp (tA)$ for each $t\ge 0$.
If $Re (z):=(z+{\tilde z})/2 >|A|$, then $|\exp (-t(zI-A))|\le
\exp (t(|A|- Re (z))\to 0$ while $t\to \infty $.
For such $z\in \bf H$ due to the Lebesgue theorem:
$(zI-A)\int_0^{\infty }\exp (-t(zI-A))dt=I$ and by Lemma 2.3 there exists
$R(z;A)=\int_0^{\infty }\exp (-t(zI-A))dt$.
For each $\epsilon >0$ let $A_{\epsilon }x:=
(T(\epsilon )x-x)/\epsilon $, where $x\in X$, for which there exists
$\lim_{0<\epsilon \to 0}A_{\epsilon }x$, a set of such $x$
we denote ${\cal D}(A)$. Evidently, that ${\cal D}(A)$ is
the $\bf H$-vector space in $X$. We take on it the infinitesimal
QLO $Ax:=\lim_{0<\epsilon \to 0}A_{\epsilon }x$.
Considering  $\bf H$ as BS over $\bf R$ we get analogs
of Lemmas 3,4,7, Corollaries 5, 9 and Theorem 10 from \S VIII.1
\cite{danschw}, moreover, ${\cal D}(A)$ is dense in $X$, $A$ is closed
QLO on ${\cal D}(A)$. Let $w_0:= \lim_{t\to \infty } ln (|T(t)|)/t$
and $z\in \bf H$ with $Re (z)>w_0$. For each $w_0<\delta <Re(z)$ due to
Corollary VIII.1.5 \cite{danschw} there exists a constant $M>0$
such that $|T(t)|\le M\exp (\delta t)$ for each $t\ge 0$. Then
there exists $R(z)x:=\int_0^{\infty }\exp (-t(zI-A))xdt$ for each
$x\in X$ and $Re (z)>w_0$, consequently, $R(z)x\in {\cal D}(A)$.
Let $T_z$ be QLO corresponding to $z^{-1}A$ instead of $T$ for $A$,
where $0\ne z\in \bf H$, moreover, ${\cal D}(A)={\cal D}(z^{-1}A)$.
Then $z^{-1}A\int_0^{\infty }\exp (-t(I-z^{-1}A))xdt=
\int_0^{\infty }\exp (-t(I-z^{-1}A)z^{-1}Axdt$, consequently,
$R(z)(zI-A)x=x$ for each $x\in {\cal D}(A)$ and $R(z)=R(z;A)$.
Therefore, $R(z;A)x=\int_0^{\infty }\exp (-t(zI-A))xdt$ for each
$z\in \rho (A)$ and $x\in X$.
\par With the help of $2.7.(iii)$ for QLO $A$ there exists QLO $B$
such that $A={\bf i}B$, where $B={ {B_i\quad {B_k-{\bf i}B_j}} \choose
{{B_k+{\bf i}B_j}\quad -B_i} }$. Since $U(t)U(t)^*=U(t)^*U(t)=I$,
then $A$ commutes with $A^*$ and $\exp (t(A+A^*))=I$.
From $R(z;B)^*=R({\tilde z},B)$ it follows, that $B=B^*$.
if ${\hat E}$ it is the decomopsition of the identity for $B$
and $V(t):=\exp ({\bf i}tB)$,
by Theorem 2.30 $<V(t)x;y>=\int^{\infty }_{-\infty }
<{\hat E}(dz).\exp ({\bf i}tz)x;y>$, then due the Fubini theorem
$\int_0^{\infty }<V(t).\exp (-bt)x;y>dt=$ $\int_0^{\infty }
\int^{\infty }_{-\infty }<{\hat E}(dz).\exp (-(b-{\bf i}z)t)x;y>dt$
$=\int^{\infty }_{-\infty }<{\hat E}(dz).(b-{\bf i}z)^{-1}x;y>$
$=R(b;{\bf i}B)x;y>$ for each $b\in \bf H$ with $Re (b)>0$.
Therefore, $\int_0^{\infty }<V(t).\exp (-bt)x;y>dt=
\int_0^{\infty }<U(t).\exp (-bt)x;y>dt$ while $Re (b)>0$.
In view of Lemma $VIII.1.15$ $<V(t).\exp (-\epsilon t)x;y>=
<U(t).\exp (-\epsilon t)x;y>$ for each $t\ge 0$ and $Re(b)>0$,
consequently, $U(t)=V(t)$.
\par {\bf 2.34. Notations.} Let $X$ be a $\bf H$-linear locally convex space.
Consider left, right and two-sided $\bf H$-linear spans of a family of
vectors $\{ v^a: a\in {\bf A} \} $,
where $span_{\bf H}^l \{ v^a: a\in {\bf A} \} :=
\{ z\in X: z=\sum_{q_a\in \bf H; a\in \bf A} q_av^a \} $;
$span_{\bf H}^r \{ v^a: a\in {\bf A} \} :=
\{ z\in X: z=\sum_{q_a\in \bf H; a\in \bf A} v^aq_a \} $;
$span_{\bf H} \{ v^a: a\in {\bf A} \} :=
\{ z\in X: z=\sum_{q_a, r_a \in \bf H; a\in \bf A} q_av^ar_a \} $.
\par {\bf 2.35. Lemma.} {\it In the notation of \S 2.34
$span_{\bf H}^l \{ v^a: a\in {\bf A} \} =
span_{\bf H}^r \{ v^a: a\in {\bf A} \} =
span_{\bf H} \{ v^a: a\in {\bf A} \} .$ }
\par {\bf Proof.}  In view of continuity of the additivity and
multiplication of vectors in $X$ and using convergence of
a net of vectors it is sufficient to prove the statement of the lemma
for a finite set $\bf A$.
Then the space $Y:=span_{\bf H} \{ v^a: a\in {\bf A} \} $
is finite-dimensional over $\bf H$ and evidently, that left and right
$\bf H$-linear spans are contained in it. Then in $Y$ it can be
chosen a basis over $\bf H$
and each vector can be written in the form $v^a= \{ v^a_1,...,v^a_n \} $,
where $n\in \bf N$, $v^a_s\in \bf H$. Each quaternion $q\in \bf H$
can be written in the form of $4\times 4$ real matrix,
therefore, for each vector $y\in Y$ there exist matrices $A$ and $B$,
elements of which belong to $\bf H$ such that $AV=y$ and $WB=y$,
where $W=\{ v^a: a\in {\bf A} \} $, $V=W^t$ is the transposed matrix,
since $A$, $B$, $W$ and $V$ can be written in the block form
over $\bf R$ with blocks $4\times 4$. Therefore, $span_{\bf H}^l
\{ v^a: a\in {\bf A} \} \cap span_{\bf H}^r \{ v^a: a\in {\bf A} \}
\supset span_{\bf H} \{ v^a: a\in {\bf A} \} $, that together with
the inclusion $span_{\bf H}^r \{ v^a: a\in {\bf A} \}
\cup span_{\bf H}^l \{ v^a: a\in {\bf A} \} \subset
span_{\bf H} \{ v^a: a\in {\bf A} \} $ proved above
leads to the statement of this lemma.
\par {\bf 2.36. Lemma.} {\it Let $X$ be HS over $\bf H$,
$\bf X_{\bf R}$ be the same space considered over the field
$\bf R$. A vector $x\in X$ is orthogonal to a $\bf H$-linear
subspace $Y$ in $X$ relative to the $\bf H$-valued scalar
product in $X$ if and only if $x$ is orthogonal to
$Y_{\bf R}$ relative to the scalar rpoduct in $X_{\bf R}$.
The space $X$ is isomorphic to the standard HS $l_2(\alpha ,{\bf H})$
over $\bf H$ of converging by the norm sequences $v=\{ v^a:
a\in \alpha \} $ with the scalar product $<v;w>:=
\sum_a{\tilde v}^aw_a$, moreover, $card (\alpha )\aleph _0=w(X)$, where
$card (\alpha )$ is the cardinality of the set $\alpha $,
$\aleph _0=card ({\bf N})$.}
\par {\bf Proof.} In view of Lemma 2.35 and transfinite induction
in $Y$ there exists a $\bf H$-linearly independent system
of vectors $\{ v^a: a\in {\bf A } \} $ such that
$span_{\bf H}^r \{ v^a: a\in {\bf A} \} $ is everywhere dense in $Y$.
In another words in $Y$ there exists a Hamel basis over $\bf H$.
A vector $x$ is by definition orthogonal to $Y$ if and only if
$<v;x>=0$ for each $v\in Y$, that is equivalent to
$<v^a;x>=0$ for each $a\in \bf A$. The space $X_{\bf R}$
is isomorphic with the direct sum $X_e\oplus iX_i\oplus jX_j\oplus kX_k$,
where $X_e$, $X_i$, $X_j$ and $X_k$ are pairwise isomorphic
HS over $\bf R$. The scalar product $<x;y>$ in $X$ can be written in
the form
\par $(i)$ $<x;y>=\sum_{m,n\in \{ e,i,j,k \} } <x_m;y_n>{\tilde m}n$, \\
where $<x_m;y_n>\in \bf R$ due to $2.16.(3,5)$. Then
the scalar product $<x;y>$ in $X$ induces the scalar product
\par $(ii)$ $<x;y>_{\bf R}:=\sum_{m\in \{ e,i,j,k \} } <x_m;y_m>$ \\
in $X_{\bf R}$. Therefore, from the orthogonality of
of $x$ to the subspace $Y$ relative to $<x;y>$ it follows orthogonality
of $x$ to the subspace $Y_{\bf R}$ relative to $<x;y>_{\bf R}$.
In view of Lemma 2.35 from $y\in Y$
it follows, that $my_m\in Y$ for each $m\in \{ e,i,j,k \} $.
Then from $<x;y_m>_{\bf R}=0$ for each $y\in Y$ and $m$
due to 2.16.(5) it follows $<x;y>=0$ for each $y\in Y$.
Then by theorem about transfinite induction \cite{eng}
in $X$ there exists an orthonormal basis over $\bf H$,
in which each vector can be represented in the form of a converging
series of left (or right) $\bf H$-linear combinations of basis vectors.
For each $x\in X$ in view of normability of $X$
the base of neighbourhoods is countable and for the topological density
we have the equality $d(X)=card (\alpha )\aleph _0$, since
$\bf H$ is separable, hence $w(X)=d(X)$.
From this the last statement of this lemma follows.
\par {\bf 2.37. Lemma.} {\it For each QLO $T$ in HS $X$ over $\bf H$
an adjoint operator $T^*$ in $X$ relative to a $\bf H$-scalar
product coinsides with an adjoint operator $T^*_{\bf R}$
in $X_{\bf R}$ relative to a $\bf R$-valued scalar product in $X_{\bf R}$.}
\par {\bf Proof.} Let ${\cal D}(T)$ be a domain of operator
$T$, which is dense in $X$. In view of Formulas $2.36.(i,ii)$
and the existense of the automorphisms $z\mapsto zm$ in $\bf H$
for each $m\in \{ e,i,j,k \} $ it follows that the continuities of
$<Tx;y>$ and $<Tx;y>_{\bf R}$ by $x\in {\cal D}(T)$ are equivalent,
therefore, in view of Lemma 2.35 the family of all $y\in X$,
for which $<Tx;y>$ is continuous by $x\in {\cal D}(T)$
forms a $\bf H$-linear subspace in $X$ and it is the same
relative to $<Tx;y>_{\bf R}$, that is the domain ${\cal D}(T^*)$ of the
operator $T^*$. Then the adjoint operator $T^*$ is defined
by the equality $<Tx;y>=:<x;T^*y>$, while $T^*_{\bf R}$ is given by
$<Tx;y>_{\bf R}=<x;T^*_{\bf R}y>_{\bf R}$,
where $x\in {\cal D}(T)$ and $y\in {\cal D}(T^*)$.
In view of Formulas $2.36.(i,ii)$ $<x_m;(T^*y)_m>=<x_m;(T^*y)_m>_{\bf R}$
for each $x\in {\cal D}(T)$, $y\in {\cal D}(T^*)$
and $m\in \{ e,i,j,k \} $.
Since due to Proposition 2.20 and Lemma 2.35
${\cal D}(T)$ and ${\cal D}(T^*)$ are $\bf H$-linear,
then automorphisms of the field $\bf H$ given above lead to
$T^*=T^*_{\bf R}$.
\par {\bf 2.38. Definition.} A bounded QLO $P$ in HS $X$ over
$\bf H$ is called a partial $\bf R$- (or $\bf H$-) isometry,
if there exists a closed $\bf R$- (or $\bf H$-) linear subspace
$Y$ such that $\| Px \| = \| x \| $ for each $x\in Y$ and
$P(Y^{\perp }_{\bf R})= \{ 0 \} $ (or $P(Y^{\perp })= \{ 0 \} $)
respectively,
where $Y^{\perp }:=\{ z\in X: \quad <z;y>=0 \quad \forall y\in Y \} $,
$Y^{\perp }_{\bf R}:=\{ z\in X_{\bf R}: \quad <z;y>_{\bf R}=0
\quad \forall y\in Y \} $.
\par {\bf 2.39. Theorem.} {\it If $T$ is closed QLO
in HS $X$ over $\bf H$, then $T=PA$, where $P$ is a partial
$\bf R$-isometry on $X_{\bf R}$ with a domain $cl (Range (T^*))$
and $A$ is a self-adjoint QLO such that $cl (Range (A))=
cl (Range (T^*))$. If $T$ is $\bf H$-linear, then $P$ is a partial
$\bf H$-isometry.}
\par {\bf Proof.} In view of the spectral Theorem 2.28 a self-adjoint
QLO $T$ is positive if and only if its spectrum
$\sigma (T)\subset [0, \infty )$
(see also Lemma XII.7.2 \cite{danschw}).
In the field $\bf H$ each polynomial has a root (see Theorem 3.17
\cite{luoyst}). Therefore, $T$ is a positive self-adjoint QLO,
then there exists a unique positive QLO $A$ such that $A^2=T$
(see also Lemma XII.7.2 \cite{danschw}).
Then there exists the positive square root $A$ of the operator
$T^*T$. Moreover, $A$ is $\bf H$-linear, if $T$ is
$\bf H$-linear. Put $SAx=Tx$ for each $x\in {\cal D}(T^*T)$
and $V$ be an isometrical extension of $S$ on $cl (Range (A))$.
The space $cl (Range (A))$ is $\bf R$-linear.
If $A$ is in addition left- (or right-) $\bf H$-linear,
then $cl (Range (A))$ is a $\bf H$-linear subspace due to Lemma
2.35. In view of Lemma 2.36 there exists the perpendicular projector
$E$ from $X$ on $cl (Range (A))$,
moreover, $E$ is $\bf H$-linear, if $cl (Range (A))$ is the
$\bf H$-linear subspace. Then put $P=VE$.
From $<Ax;Ax>=<Tx;Tx>$ for each $x\in {\cal D}(T^*T)$ it follows, that
$PAx=Tx$ for each $x\in {\cal D}(T^*T)$.
The rest of the proof can be done analogously to the proof
of Theorem $XII.7.7$ \cite{danschw} with the help of
Lemmas 2.35-37.
\par {\bf 2.40. Note and Definition.} Apart from the case of
$\bf C$ nontrivial polynomials of quaternion variables can have
roots, which are not points< but closed submanifolds in
$\bf H$ with dimensions over $\bf R$ from $0$ up to $3$ (see \cite{luoyst}).
\par A closed subset $\lambda \subset \sigma (T)$ is called
an isolated subset of the spectrum, if there exists a neighbourhood
$U$ of a subset $\lambda $ such that $\sigma (T)\cap U=\lambda $.
An isolated subset $\lambda $ of a spectrum $\sigma (T)$
is called a pole of a spectrum (of order $p$), if $R(z;T)$ has
zero on $\lambda $ (of order $p$, that is, each $z\in \lambda $
is zero of order $0<p(z)\le p$ for $R(z;T)$ and
$\max_{z\in \lambda }p(z)=p$). A subset $\lambda $ which
is clopen (closed and open simultaneously) in $\sigma (T)$
is called a spectral set.
Let $\eta _1$ be a closed rectifiable path in $U$ encompassing
$\lambda $ and not intersecting with $\lambda $, characterized
by a vector $M_1\in \bf H$, $|M_1|=1$, $M_1+{\tilde M}_1=0$
(see Theorem 3.22 \cite{luoyst}), denote
$$(i)\quad \phi _n(z,T):=(2\pi )^{-1}
\{ \int_{\eta _1} R(\zeta ;T)((\zeta -a)^{-1}(z-a))^n(\zeta -a)^{-1} d\zeta
)M_1^{-1} \} ,$$
where $\eta _1\subset B({\bf H},a,R)\setminus B({\bf H},a,r)$,
$B({\bf H},a,R)\subset U$, $\lambda \subset B({\bf H},a,r)$,
$0<r<R<\infty $.
We say, that an index of $\lambda $ is equal to $p$ if and only if
there exists a vector $x\in X$ such that
$$(ii)\quad (zI-T)^{s_1}v_1{\hat E}(\delta (z);T)v_2...(zI-T)^{s_m}v_{2m-1}
{\hat E}(\delta (z);T)x=0$$
for each $z\in \lambda $ and each $0\le s_n\in \bf Z$ with
$s_1+...+s_m=p$ and each $v_1,..,v_{2m-1}$,
where $v_1=v_1(\delta ,T)\in \bf H$,...,$v_{2m-1}=v_{2m-1}(\delta ,T)
\in \bf H$, $m\in \bf N$, $\delta :=\delta (z)\ni z$,
$\delta (z)\in {\cal B}(\lambda )$,
while expression in $(ii)$ is not equal to zero for some $z\in \lambda $
and $s_1,...,s_m$ with $s_1+...+s_m=p-1$.
\par {\bf 2.41. Theorem.} {\it A subset $\lambda $ is a pole
of order $p$ of QLO $T\in L_q(X)$ for $U=B({\bf H},\alpha ,R')$,
$0<R<R'<\infty $ in Definition 2.40, where $0<r<\infty $ if and only if
$\lambda $ has an index $p$.}
\par {\bf Proof.} Choose with the help of a homotopy relative to
$U\setminus \lambda $ closed paths $\eta _1$ and $\eta _2$
homotopic to $\gamma _1$ and $\gamma _2$, moreover, $\inf _{\theta }
|\eta _1(\theta )| > \sup_{\theta }|\eta _2(\theta )|$,
where $\gamma _1$ and $\gamma _2$ are chosen as in Theorem
$3.22$ \cite{luoyst} (see also Theorem $3.9$ there), $\theta \in [0,1]$.
In view of Theorem 3.22 \cite{luoyst} the quaternion Loran decomposition
of $R(z;T)$ in the neighbourhood $B({\bf H},a,R)\setminus B({\bf H},a,r)$
has the form $R(z;T)=\sum_{n=0}^{\infty }(\phi _n(z,T)+ \psi _n(z,T))$,
where $\phi _n$ is given by Formula $2.40.(i)$ and
$$(i)\quad \psi _n(z,T):=(2\pi )^{-1}
\{ \int_{\eta _2} R(\zeta ;T)(z-a)^{-1}((\zeta -a)(z-a)^{-1})^n
d\zeta )M_2^{-1} \} .$$ If $\lambda $ is a pole of order
$p$, then $\phi _p=0$ and $\phi _{p-1}\ne 0$, then there exists
$x\in X$ such that
$$(ii)\quad \phi _{p(z)}(z,T)x=0 \mbox{ for each }z\in \lambda ,
\mbox{ and }$$
$$(iii)\quad \phi _{p-1}(z,T)x\ne 0\mbox{ for some }z\in \lambda .$$
An analogous decompositions with the corresponding $\phi _n$
are true for the products $f(T)R(z;T)g(T)$,
where $f$ and $g$ are quaternion holomorphic functions on a neighbourhood
$\sigma (T)$ not equal to zero everywhere on $\lambda $.
Functions $\phi _n$ for $R(z;T)$
can be approximated with any precision in the strong operator topology
in the form of left $\bf H$-linear combinations of functions
taking part in $2.40.(ii)$ due to Lemma 2.35 and the definition
of the quaternion line integral along a rectifiable path, since
while $|\xi |> \sup |\chi |$ the series for
$R(\xi ;T_{\chi })$ converges in the uniform operator topology
for each spectral set $\chi $ of the spectrum $\sigma (T)$,
where $T_{\chi }=T|_{X_{\chi }}$, $X_{\chi }:={\hat E}_e(\chi ;T)X$.
The variation of $f$ and $g$ implies, that the index of $\lambda $
is not less than $p$.
Vice versa, let Conditions $(ii)$ be satisfied for some $n$.
The resolvent $R(z;T)x$ is regular on ${\bf H}\setminus B({\bf H},a,r)$
and
$$x=(2\pi )^{-1}\{ \int_{\eta }R(\zeta ;T)x d\zeta \} M^{-1}=
(2\pi )^{-1}\{ \int_{\eta _2}R(\zeta ;T)x d\zeta \} M_2^{-1}=
\omega (T)x,$$
where $\omega (T)$ is the function equal to $1$ on a neighbourhood of
$\lambda $ and equal to zero on ${\bf H}\setminus U$,
$\eta $ is the corresponding closed rectifiable path encompassing
$\sigma (T)$ and characterized by $M\in \bf H$, $|M|=1$,
$M+{\tilde M}=0$. Then due to $2.40.(ii)$ $\phi _{p(z)}(z,T)x=0$
for each $z\in \lambda $.
\par {\bf 2.42. Note.} An isolated point $\lambda $
of a spectrum $\sigma (T)$ for normal QLO
$T\in L_q(X)$ on HS $X$ over $\bf H$
may not have eigenvectors because of non-commutativity of
a projection-valued measure $\hat E$ apart from the case
of linear operators on HS over $\bf C$.

\par {\underline {Acknowledgement.}} The author is sincerely
grateful to Prof. H. de Groote for his interest to this work,
discussions and hospitality at Department of Applied Mathematics
of Frankfurt-am-Main University in Autumn 2002 - January 2003.
\par Address: Mathematical Department, Brussels University, 
V.U.B., Pleinlaan 2, Brussels 1050, Belgium
\end{document}